\documentclass[a4paper,12pt,oneside]{amsart}
\usepackage{graphicx}
\usepackage{booktabs}
\usepackage[utf8]{inputenc}
\usepackage[T1]{fontenc}
\usepackage{booktabs}
\newcommand{\functionname}[1]{%
  \par\medskip\noindent\codebreak{#1}\par\smallskip
}

\newcommand{\codebreak}[1]{\texttt{\seqsplit{#1}}}
\usepackage{seqsplit}

\usepackage[
  a4paper,
  left=3cm,
  right=3cm,
  top=4cm,
  bottom=4cm,
  twoside=false 
]{geometry}

\usepackage{amsmath,amsthm,verbatim,amssymb,amsfonts,amscd,graphicx,mathrsfs,array,float,mathtools,tikz-cd,enumitem}
\usepackage{tabularx}
\usepackage{diagbox}

\numberwithin{equation}{section}

\newtheorem{theorem}{Theorem}[section]
\newtheorem{corollary}[theorem]{Corollary}

\newtheorem{conjecture}[theorem]{Conjecture}

\theoremstyle{definition}
\newtheorem{definition}[theorem]{Definition}

\theoremstyle{remark}
\newtheorem{remark}[theorem]{Remark}
\newtheorem{example}[theorem]{Example}

\usepackage{thmtools}
\usepackage{thm-restate}
\usepackage{hyperref}
\usepackage{cleveref}

\newcommand{\PP}{\mathbb{P}}

\newcommand{\CC}{\mathbb{C}}

\DeclareMathOperator{\Sym}{Sym}

\title{Families of Smooth Fano Fourfolds of Picard Rank 1 without Bott Vanishing}

\author[J. Wang]{Jiahe Wang}
\address{UCLA Mathematics Department, Box 951555, Los Angeles, CA 90095-1555, USA
}
\email{jiahewang@math.ucla.edu}
\begin{document}
\pagestyle{plain}

\begin{abstract}
We show that \(\chi(X,T_X)<0\) for the currently known
families of smooth Fano fourfolds of Picard rank \(1\) and index \(1\). Combining this with the known Picard rank $1$ index $> 1$ cases, we show that
among all currently known smooth Fano fourfolds of Picard rank $1$, the only
variety satisfying Bott vanishing is the projective space. By a result of Kawakami--Totaro, the existence of an endomorphism of degree greater than 1 implies Bott vanishing. Therefore, among the currently known
smooth Fano fourfolds of Picard rank \(1\), any variety admitting an
endomorphism of degree greater than 1 must be \(\mathbb P^4\). Together with Burt Totaro, we develop new Schubert2 functions for symmetric and skew-symmetric degeneracy loci, and weighted projective spaces. 
\end{abstract}

\maketitle
\section{Introduction}
The following conjecture is well-known.

\begin{conjecture}
Let \(X\) be a smooth Fano variety of Picard number \(\rho(X)=1\) over an
algebraically closed field of characteristic zero. If \(X\) admits an
endomorphism of degree greater than \(1\), then \(X\simeq \mathbb P^n\).
\end{conjecture}
The conjecture is known in several cases: when \(\dim X=3\) \cite{amerik,arv,hwang}, when \(\dim X=4\) and \(X\) has Fano index greater than \(1\) \cite{shao}, when \(X\) is a hypersurface \cite{paranjape,beauville}, and when \(X\) is a homogeneous space \cite{paranjape}. Kawakami and Totaro \cite{totaro} prove a general implication from int-amplified endomorphisms to Bott vanishing. Int-amplified is a stronger condition than an endomorphism of degree greater than 1; when \(\rho(X)=1\), they are equivalent. As applications, they prove the above cases in arbitrary characteristic.

\begin{theorem}[{\cite[Theorem C]{totaro}}]
Let \(X\) be a normal projective variety over a perfect field \(k\). Suppose
that \(X\) admits an int-amplified endomorphism whose degree is invertible in
\(k\). Then \(X\) satisfies Bott vanishing for ample Weil divisors. That is,
\[
H^i\bigl(X,\Omega_X^{[j]}(A)\bigr)=0
\]
for every \(i>0\), \(j\geq 0\), and \(A\) an ample Weil divisor.
\end{theorem}
In this paper we apply Kawakami--Totaro's Bott vanishing implication to 32 families of known smooth Fano fourfolds with Picard rank 1 and index 1 over $\CC$. The
families considered here are:
\begin{itemize}
    \item \(13\) families realized as complete intersections in weighted
    projective spaces \cite{kuchler,Przyjalkowski};
    \item \(13\) families realized as zero loci of homogeneous vector bundles on
    Grassmannians \cite{kuchle};
    \item \(6\) additional families realized as skew-symmetric degeneracy loci
    \cite{qureshi,kuchler}.
\end{itemize}

We show that, for each case, the Euler characteristic $\chi(X, T_X)$ is negative and therefore disprove Bott vanishing for the families: indeed, if Bott vanishing held, then since \(T_X \cong \Omega_X^{3}\otimes (-K_X)\) and \(-K_X\) is ample, we would have \(H^i(X,T_X)=0\) for all \(i>0\), and hence \(\chi(X,T_X)=h^0(X,T_X)\ge 0\). By Kawakami--Totaro, this shows that each variety in these families has no endomorphism of degree greater than 1. Together with the 16 Picard rank 1 index greater than 1 families \cite[Proposition 3.12]{totaro}, this verifies the endomorphism conjecture for the currently known smooth Fano fourfolds of Picard rank \(1\). \begin{theorem} For each of the 32 currently known families of smooth Fano fourfolds of Picard rank 1 and index 1, one has $\chi(X,T_X)<0$.  In particular, none of these varieties satisfies Bott vanishing. \end{theorem} \begin{corollary} Among the known smooth Fano fourfolds of Picard rank 1, any variety admitting an endomorphism of degree greater than 1 is isomorphic to $\PP^4$. \end{corollary} 

Together with Burt Totaro, we develop several computational functions for the
Macaulay2 \cite{M2} package \codebreak{Schubert2} \cite{Schubert2Source}, and
apply them to known families of Fano fourfolds. The new functions include
\codebreak{skewDegeneracyLocus}, \codebreak{skewDegeneracyLocus2}, and
\codebreak{skewKernelBundle} for skew-symmetric degeneracy loci, together with the
analogous functions for symmetric degeneracy loci. These functions make it possible to compute
characteristic classes and Euler characteristics for degeneracy loci. We also include
\codebreak{weightedProjectiveBundle} and \codebreak{weightedProjectiveSpace}, which
allow intersection-theoretic computations in
weighted projective spaces.

\subsection*{Acknowledgments}

I am grateful to Burt Totaro for proposing this topic and for helpful guidance and comments. I am grateful to Enrico Fatighenti for helpful comments on
families of Fano fourfolds.

\section{Preliminaries}
\begin{definition}
Let \(P\) be a variety, let \(E\) be a vector bundle on \(P\), and let \(L\) be a line bundle on \(P\). A morphism
\(\varphi:E\to E^\vee\otimes L\) is called \emph{skew-symmetric} if the corresponding section of \(E^\vee\otimes E^\vee\otimes L\) lies in \(\wedge^2E^\vee\otimes L\). Similarly, \(\varphi\) is called \emph{symmetric} if the corresponding section lies in \(\operatorname{Sym}^2E^\vee\otimes L\).
\end{definition}

\begin{definition}
Let \(\varphi:E\to E^\vee\otimes L\) be a symmetric or skew-symmetric morphism. For an integer \(k\), the \emph{rank degeneracy locus} associated to \(\varphi\) is
\[
D_k(\varphi):=\{x\in P:\operatorname{rank}(\varphi_x)\leq k\}.
\]
Locally, after trivializing \(E\) and \(L\), this is the locus where the corresponding symmetric or skew-symmetric matrix has rank at most \(k\). In the symmetric case, \(D_k(\varphi)\) is cut out by the \((k+1)\times(k+1)\) minors of this symmetric matrix. In the skew-symmetric case, the rank is always even, and \(D_{2s}(\varphi)\) is locally cut out by the Pfaffians of size \(2s+2\). The degeneracy loci of skew-symmetric morphisms are also called \emph{Pfaffian subvarieties}.
\end{definition}

\begin{remark}
Let \(V\) be an \(n\)-dimensional vector space, and let
\[
\varphi:V\otimes\mathcal O_{\mathbb P(\wedge^2V^\vee)}
\to V^\vee\otimes\mathcal O_{\mathbb P(\wedge^2V^\vee)}(1)
\]
be the universal skew-symmetric morphism. In this special case, one writes
\(
\operatorname{Pf}(2s,n):=D_{2s}(\varphi)\subset \mathbb P(\wedge^2V^\vee).
\)
In particular, \(\operatorname{Pf}(2,n)\) is the Grassmannian \(\operatorname{Gr}(2,V^\vee)\).
\end{remark}
\begin{remark}
To compute algebraic invariants of a degeneracy locus $X$, it is often enough to work on an open smooth neighborhood of \(X\). Thus the ambient space $P$ in this paper may have singularities away from \(X\).
\end{remark}
\section{New Schubert2 functions}

Together with Burt Totaro, we develop several functions to be added for the existing Macaulay2 package \texttt{Schubert2}. These functions are designed to compute
characteristic classes and Euler characteristics of varieties defined by
symmetric or skew-symmetric degeneracy conditions, as well as complete
intersections in weighted projective spaces.
\subsection*{Method}
The Macaulay2 code of the functions was initially generated with the assistance of ChatGPT and was subsequently revised by the
author. The Macaulay2 descriptions of the functions and their examples were
written by Burt Totaro.

\functionname{skewDegeneracyLocus}

Let \(\varphi:E\to E^\vee\otimes L\) be a skew-symmetric map over $P$ with $E$ of rank $n$. For any even number $k$, this function constructs the resolution of the rank $k$
degeneracy locus of \(\varphi\) as a zero locus \(Z\) of a section of a vector bundle on the relative Grassmannian
\(G=\operatorname{Gr}(n-k,E)\) over $P$. If
\(0\to S\to p^*E\to Q\to0\) is the tautological sequence on \(G\), then the
bundle defining \(Z\) has a natural filtration whose associated graded bundle is
\(
F_{\mathrm{skew}}:=
(S^\vee\otimes Q^\vee\otimes p^*L)\oplus
(\wedge^2S^\vee\otimes p^*L).
\)
Thus \(F_{\mathrm{skew}}\) may be used for Chern-class and \(K\)-theory
computations.

This formula follows from Example 14.4.11 in Fulton \cite{fulton}. The rank condition \(\operatorname{rank}\varphi_x\leq k\) is equivalent to
\(\dim\ker(\varphi_x)\geq n-k\). Thus one considers the relative Grassmannian
\(p:G=\operatorname{Gr}(n-k,E)\to P\). A point of \(G\) is a pair
\((x,S_x)\), where \(x\in P\) and \(S_x\subset E_x\) is an
\((n-k)\)-dimensional subspace. Let \(Z\subset G\) be the locus of pairs
\((x,S_x)\) such that \(S_x\subset\ker(\varphi_x)\). Then \(Z\) is the
vanishing locus of the composition
\(S\to p^*E\xrightarrow{p^*\varphi}p^*E^\vee\otimes p^*L\), equivalently the vanishing of a section of \(S^\vee\otimes p^*E^\vee\otimes p^*L\). The tautological sequence gives a filtration of
\(S^\vee\otimes p^*E^\vee\otimes p^*L\), and so \(S^\vee\otimes p^*E^\vee\otimes p^*L\) may be represented, for Chern-class computations, by
\(F_{\mathrm{skew}}=(S^\vee\otimes Q^\vee\otimes p^*L)\oplus
(\wedge^2S^\vee\otimes p^*L)\). 

If the degeneracy locus \(D_k(\varphi)\) is smooth,
contains no lower-rank stratum, and has expected codimension
\(\binom{n-k}{2}\) in $P$, then \(Z\to D_k(\varphi)\) is an isomorphism. 
Moreover, the induced section of the vector bundle is regular, so
\([Z]=c_{\operatorname{top}}(F_{\mathrm{skew}})\), and also the tangent-normal sequence gives
\(
[T_Z]=[T_G|_Z]-[F_\mathrm{skew}|_Z]
\)
in \(K(Z)\).

\functionname{symmetricDegeneracyLocus}

Similar to the skew-symmetric case, this function constructs the resolution of a symmetric
degeneracy locus as the zero locus $Z$ of a section of a vector bundle on a relative Grassmannian. For an integer \(k\), the
function constructs a zero locus \(Z\) on \(G=\operatorname{Gr}(n-k,E)\). If
\(0\to S\to p^*E\to Q\to0\) is the tautological sequence on \(G\), then the function outputs the expected class of \(Z\), \([Z]=c_{\operatorname{top}}(F_{\mathrm{sym}})\), where
\(F_{\mathrm{sym}}=(S^\vee\otimes Q^\vee\otimes p^*L)\oplus
(\operatorname{Sym}^2S^\vee\otimes p^*L)\). 

\functionname{skewDegeneracyLocus2}

This function computes the class of a skew-symmetric degeneracy locus 
on the original variety \(P\). Let \(\varphi:E\to E^\vee\otimes L\) be a general
skew-symmetric map on a smooth variety $P$, and let \(n=\operatorname{rank}E\) and \(d=n-k\). For even
\(k\), the expected codimension of \(D_k(\varphi)\) is \(d(d-1)/2\). The
function returns the class
\([D_k(\varphi)]=s_{(d-1,d-2,\ldots,1)}(E^\vee\otimes L^{1/2})\) in
\(A^{d(d-1)/2}(P)\), where \(s_\lambda\) denotes the Schur polynomial and
\(L^{1/2}\) is understood formally. This is the skew-symmetric degeneracy
locus formula of Harris--Tu \cite[Theorems 8 and 10]{harris}; see also
J\'ozefiak--Lascoux--Pragacz \cite{jozefiak}.

\functionname{symmetricDegeneracyLocus2}

This function computes the class of a symmetric degeneracy locus directly on
the original variety \(P\). Let \(\varphi:E\to E^\vee\otimes L\) be a general symmetric map on a smooth variety $P$,
let \(n=\operatorname{rank}E\), and put \(d=n-k\). The expected codimension of
\(D_k(\varphi)\) is \(d(d+1)/2\). The function returns the class
\([D_k(\varphi)]=2^d s_{(d,d-1,\ldots,1)}(E^\vee\otimes L^{1/2})\) in
\(A^{d(d+1)/2}(P)\), where \(s_\lambda\) denotes the Schur polynomial and
\(L^{1/2}\) is understood formally. This is the symmetric degeneracy locus
formula of Harris--Tu \cite[Theorems 1 and 10]{harris}; see also
J\'ozefiak--Lascoux--Pragacz \cite{jozefiak}.

\functionname{skewKernelBundle}

This function returns the tautological kernel bundle on the
relative-Grassmannian resolution of a skew-symmetric degeneracy locus. In the
notation above, let \(Z\subset G=\operatorname{Gr}(n-k,E)\) be the zero locus
constructed by \codebreak{skewDegeneracyLocus}, and let \(i:Z\hookrightarrow G\)
be the inclusion. Then the function returns \(K=i^*S\), where \(S\) is the
tautological subbundle on \(G\). When \(Z\to D_k(\varphi)\) is an isomorphism, this
bundle is the actual kernel bundle of the skew-symmetric morphism \(\varphi\)
along the degeneracy locus.

\functionname{symmetricKernelBundle}

This function returns the tautological kernel bundle on the
relative-Grassmannian resolution of a symmetric degeneracy locus. If
\(Z\subset G=\operatorname{Gr}(n-k,E)\) is constructed by
\codebreak{symmetricDegeneracyLocus}, and \(i:Z\hookrightarrow G\) is the
inclusion, then the function returns \(K=i^*S\), where \(S\) is the tautological
subbundle on \(G\). When \(Z\to D_k(\varphi)\) is an isomorphism, this is the actual
kernel bundle of the symmetric morphism \(\varphi\) along the degeneracy locus.

\functionname{weightedProjectiveBundle}

This function constructs the trivial weighted projective bundle over a base
variety \(X\). For weights \(a_0,\ldots,a_n\), it models
\(X\times\mathbb P(a_0,\ldots,a_n)\) with projection
\(\pi:X\times\mathbb P(a_0,\ldots,a_n)\to X\). If
\(H=c_1(\mathcal O(1))\), the intersection ring is modeled as
\(A^*(X)[H]/(H^{n+1})\), and the pushforward to \(X\) is normalized by
\(\pi_*(H^n)=1/(a_0a_1\cdots a_n)\). The relative tangent bundle is represented
by \(T_\pi=\sum_{i=0}^n\mathcal O(a_i)-\mathcal O\). This allows one to perform
intersection-theoretic and Riemann--Roch computations on weighted projective
bundles in \codebreak{Schubert2}.
 
\functionname{weightedProjectiveSpace}

This function constructs an abstract weighted projective space
\(\mathbb P(a_0,\ldots,a_n)\) for intersection-theoretic computations in
\codebreak{Schubert2}. It is the special case of
\codebreak{weightedProjectiveBundle} where the base variety is a point.

\section{Computation of $\chi(X, T_X)$ for known smooth Fano Fourfolds}

\subsection{Families of smooth Fano fourfolds}

We know there are finitely many families of smooth Fano varieties in each dimension \cite{KYM}. In dimension 4, the smooth Fano varieties of index at least two have been classified; there are $35$ deformation families and $16$ of them are of Picard rank 1 \cite[Theorem 1.2]{kp}, \cite[Theorem 2]{mukai}, \cite{wilson}.  The index 1 case remains incomplete, although there are many constructions.  The toric Fano fourfolds were classified by Batyrev \cite{batyrev}.  Coates, Kasprzyk and Prince constructed Fano fourfolds as complete intersections in toric varieties \cite{Coates}.  Kalashnikov with Coates and Kasprzyk \cite{Kalashnikov} constructed quiver flag zero loci, and all the families were later interpreted as zero loci of sections of homogeneous vector bundles on homogeneous varieties by Fatighenti, Tanturri, and Tufo \cite{Fatighenti}. Earlier examples include the families constructed by K\"uchle \cite{kuchler} as complete intersections in weighted projective spaces and zero loci of homogeneous vector bundles on Grassmannians \cite{kuchle}; see a list of 637 deformation families as zero loci of homogeneous vector bundles on products of Grassmannians by Bernardara, Fatighenti, Manivel, and Tanturri \cite{FanoFourfoldsDatabase} and its further discussion on K3 type \cite{fa}; see the classification of smooth Fano fourfolds and fivefolds that are complete intersections in weighted projective spaces by Przyjalkowski and Shramov \cite{Przyjalkowski}. Qureshi \cite{qureshi} constructed families of smooth Fano fourfolds as Pfaffian subvarieties.

In our case, we study Fano fourfolds with Picard number 1 and index 1. Most of the families above are not Picard number 1; for example, a smooth well-formed Picard rank 1 Fano complete intersection
in a toric variety is a weighted complete intersection \cite{Przy}. One sees that the currently known families of smooth Fano fourfolds of Picard rank 1 and index 1 are among the following: 
\begin{itemize} 
    \item $13$ families as complete intersections in weighted projective spaces
    \cite{kuchler,Przyjalkowski}.
    \item $13$ families as zero loci of homogeneous vector bundles on Grassmannians
    \cite{kuchle}.
    \item 6 additional families as Pfaffian subvarieties \cite{qureshi,kuchler}.
\end{itemize}

\subsection{Pfaffian subvarieties}

We recall two basic constructions of Pfaffian subvarieties.
\begin{example}
Let \(V\) be an \(n\)-dimensional vector space, with \(n\) odd. The projective space \(\mathbb P(\wedge^2 V^\vee)\) parametrizes skew-symmetric forms on \(V\). On \(\mathbb P(\wedge^2 V^\vee)\) there is a universal skew-symmetric map \(f:V\otimes\mathcal O\to V^\vee\otimes\mathcal O(1)\). For an even integer \(2\leq k\leq n-3\), the Pfaffian subvariety \(\operatorname{Pf}(k,n)\subset \mathbb P(\wedge^2 V^\vee)\) is the locus where this universal skew form has rank at most \(k\). Equivalently, it is cut out by the Pfaffians of size \(k+2\). In particular for the maximal Pfaffian $\operatorname{Pf}(n-3,n)$, it is cut out by the $n$ Pfaffians of size $n-1$. For $k=2$, \(\operatorname{Pf}(2,n)\) is the Pl\"ucker-embedded Grassmannian \(\operatorname{Gr}(2,V^\vee)\). For \(k>2\), however, \(\operatorname{Pf}(k,n)\) is generally not a Grassmannian; it is a higher Pfaffian degeneracy locus.
\end{example}

\begin{example}
We recall the weighted version in Qureshi \cite{qureshi}. Let \(a=(a_1,\ldots,a_n)\) be a vector of half-integers such that \(a_i+a_j\) is a positive integer for every \(i<j\). Assign weight \(a_i+a_j\) to the skew coordinate \(x_{ij}\). Then the same Pfaffian equations define a subvariety of the weighted projective space \(\mathbb P(a_i+a_j)_{i<j}\). In the case \(k=2\), this is the weighted Grassmannian \(w\operatorname{Gr}(2,n)\), obtained as the quotient of the punctured affine cone over the Pl\"ucker-embedded Grassmannian \(\operatorname{Gr}(2,n)\) by the \(\mathbb C^*\)-action with weights \(a_i+a_j\) on the Pl\"ucker coordinates. To get the map \(f:E\to E^\vee\otimes L\), one chooses \(E=\bigoplus_i\mathcal O(e_i)\) and \(L=\mathcal O(\ell)\) so that the \((i,j)\)-entry has degree \(\ell-e_i-e_j=a_i+a_j\). The choice of the integers $e_i$ and $l$ is not unique: different choices differ by a common twist, and therefore give the same weighted skew-symmetric matrix, the same Pfaffian equations, and the same degeneracy locus.
\end{example}

Qureshi \cite{qureshi} gives six families of Picard rank \(1\), index \(1\) Fano fourfolds. Qureshi constructs these families from weighted Grassmannians $w \operatorname{Gr}(2,5)$ by taking cones and quasilinear sections on the ambient weighted projective space. Here a quasilinear section means a hypersurface
whose degree is equal to one of the weights of the ambient weighted projective
coordinates. Taking cones adds variables which do not appear in the Pfaffian matrix, and taking a quasilinear section replaces a Pl\"ucker coordinate by a weighted homogeneous form of the same degree. Therefore the resulting subvarieties remain Pfaffian subvarieties of the same type. In two families, one intersects the Pfaffian subvariety with an extra hypersurface, after these operations. 

Among Qureshi's six Picard rank 1
families, example \(\#5\) already appeared in Küchle's list \cite{kuchle} as the zero locus of
\(\mathcal O_{\operatorname{Gr}(2,5)}(3)\oplus
\mathcal O_{\operatorname{Gr}(2,5)}(1)\) on Grassmannian. We therefore list the five remaining
Qureshi cases, $\#1$, $\#2$, $\#3$, $\#4$, and $\#6$, together with Küchle's additional degree \(14\) Pfaffian
family, which is a \(7\times7\) maximal Pfaffian subvariety, denoted as $\#A$. In the table, the degrees of
equations mean the degrees of the Pfaffians, together with the degree of the
extra hypersurface section when one is present. We note that example \#1 also appears in Küchle 
\cite{kuchler}.
\newpage
\begin{center}
\scriptsize
\begin{tabular}{c|c|c|c|c|c}

\text{\#} 
& \text{Type} 
& \text{Eq Degs, wP} 
& \text{Weight Matrix} 
& \(E\) 
& \(L\)
\\ \hline

1
& \(\operatorname{Gr}(2,5)\)
& \(X_{2,3^4}\subset \mathbb P^7\)
&
\(\begin{matrix}
1&1&1&2\\
 &1&1&2\\
 &&1&2\\
 &&&2
\end{matrix}\)
& \(\mathcal O^{\oplus 4}\oplus\mathcal O(-1)\)
& \(\mathcal O(1)\)
\\[1.5em]

2
& \(\operatorname{Gr}(2,5)\)
& \(X_{3^4,4}\subset \mathbb P^7(1^7,2)\)
&
\(\begin{matrix}
1&1&1&1\\
 &2&2&2\\
 &&2&2\\
 &&&2
\end{matrix}\)
& \(\mathcal O\oplus\mathcal O(-1)^{\oplus 4}\)
& \(\mathcal O\)
\\[1.5em]

3
& \(\operatorname{Gr}(2,5)\)

& \(X_{3^2,4^3}\subset \mathbb P^7(1^6,2^2)\)
&
\(\begin{matrix}
1&1&2&2\\
 &1&2&2\\
 &&2&2\\
 &&&3
\end{matrix}\)
& \(\mathcal O^{\oplus 3}\oplus\mathcal O(-1)^{\oplus 2}\)
& \(\mathcal O(1)\)
\\[1.5em]

4
& \(\operatorname{Gr}(2,5)\)
& \(X_{4^5}\subset \mathbb P^7(1^5,2^3)\)
&
\(\begin{matrix}
2&2&2&2\\
 &2&2&2\\
 &&2&2\\
 &&&2
\end{matrix}\)
& \(\mathcal O^{\oplus 5}\)
& \(\mathcal O(2)\)
\\[1.5em]

6
& \(\operatorname{Gr}(2,5)\cap\mathcal H\)
& \(X_{2^5,(4)}\subset \mathbb P^8(1^8,2)\)
&
\(\begin{matrix}
1&1&1&1\\
 &1&1&1\\
 &&1&1\\
 &&&1
\end{matrix}\)
& \(\mathcal O^{\oplus 5}\)
& \(\mathcal O(1)\)
\\[1.5em]

A
& \(\operatorname{Pf}(4,7)\)
& \(X_{3^7}\subset \mathbb P^7\)
&
\(\begin{matrix}
1&1&1&1&1&1\\
 &1&1&1&1&1\\
 &&1&1&1&1\\
 &&&1&1&1\\
 &&&&1&1\\
 &&&&&1
\end{matrix}\)
& \(\mathcal O^{\oplus 7}\)
& \(\mathcal O(1)\)
\\

\end{tabular}
\end{center}

Qureshi gives the construction of example \#3 and \#6. The other examples may be implied by the weight matrices and their ambient weighted projective spaces.

\subsection*{Example \#1}
We use the weight \(w=\frac12(1,1,1,1,3)\) to get the
embedding
\(
\operatorname{wGr}(2,5)\hookrightarrow \mathbb P^{9}(1^6,2^4).
\)
Let \(Y_1\subset \mathbb P^{11}(1^8,2^4)\) be a projective cone over
\(\operatorname{wGr}(2,5)\) with vertex \(\mathbb P^1\), i.e. we add two variables
of weight \(1\) to the ambient weighted projective space which are not involved in
any defining equations of \(\operatorname{wGr}(2,5)\). We take a complete intersection
of \(Y_1\) with four general quasilinear sections of degree \(2\) to get a Fano
4-fold
\(
X\subset \mathbb P^7.
\)

\subsection*{Example \#2.}
We use the weight
\(
w=\frac{1}{2}(0,2,2,2,2)
\)
to get an embedding
\(
\mathrm{w}G\hookrightarrow \mathbb P^9(1^4,2^6).
\)
Let \(Y_1\subset \mathbb P^{12}(1^7,2^6)\) be a projective cone over
\(\mathrm{w}G\) with vertex \(\mathbb P^2\), i.e. we add three variables of
weight \(1\) to the ambient weighted projective space which are not involved in
any defining equations of \(\mathrm{w}G\). We take a complete intersection of
\(Y_1\) with five general quasilinear sections of degree \(2\) to get a Fano 4-fold
\(
X\subset \mathbb P^7(1^7,2).
\)

\subsection*{Example \#3({\cite[3.1]{qureshi}}). }

For \(w=\frac{1}{2}(1,1,1,3,3)\) we get the embedding
\(
\mathrm{w}G \hookrightarrow \mathbb P^9(1^3,2^6,3).
\)
Let \(Y_1\subset \mathbb P^{12}(1^6,2^6,3)\) be a projective cone over
\(\mathrm{w}G\) with vertex \(\mathbb P^2\), i.e. we add three variables of
weight 1 to the ambient \(\mathbb P^9(w_i)\) which are not involved in any
defining equations of \(\mathrm{w}G\).
We take a complete intersection of \(Y_1\) with 4 general quadrics to get a
Fano 5-fold
\(
Y_2\subset \mathbb P^8(1^6,2^2,3).
\) As a last step we take an
intersection of \(Y_2\) with a general cubic to get a Fano 4-fold
\(
X\subset \mathbb P^7(1^6,2^2).
\)

\subsection*{Example \#4.}
We use the weight
\(
w=\frac{1}{2}(2,2,2,2,2)
\)
to get an embedding
\(
\mathrm{w}G\hookrightarrow \mathbb P^9(2^{10}).
\)
Let \(Y_1\subset \mathbb P^{14}(1^5,2^{10})\) be a projective cone over
\(\mathrm{w}G\) with vertex \(\mathbb P^4\), i.e. we add five variables of
weight \(1\) to the ambient weighted projective space which are not involved in
any defining equations of \(\mathrm{w}G\). We take a complete intersection of
\(Y_1\) with seven general quasilinear sections of degree \(2\) to get a Fano
4-fold
\(
X\subset \mathbb P^7(1^5,2^3).
\)

\subsection*{Example \#6 ({\cite[3.2]{qureshi}}).}

The Grassmannian \(\operatorname{Gr}(2,5)\) has the embedding in
\(\mathbb P^9_{(x_1,\cdots,x_{10})}\). Let \(Y_1\) be the variety obtained by taking a cone of
weight 2 over it, i.e. we have the embedding
\(
Y_1=C^2\operatorname{Gr}(2,5)\hookrightarrow \mathbb P^{10}(1^{10},2)
\)
where the new variable \(y\) of weight 2 is not involved in any defining
equations. Now
we take a general quartic section
\(
Q_4=y^2+f_4(x_i,y),\quad 1\leq i\leq 10
\)
of \(Y_1\) to get a 6-fold
\(
Y_2\subset \mathbb P^{10}(1^{10},2).
\)
We take two hyperplane sections of \(Y_2\) to get
a smooth Fano 4-fold \(X\subset \mathbb P^8(1^8,2)\).

\subsection*{Example \#A}
We use the ordinary Pfaffian format, with all skew-matrix entries of weight \(1\), to
get an embedding
\(
\operatorname{Pf}(4,7)\hookrightarrow \mathbb P(\wedge^2\mathbb C^7)=\mathbb P^{20}.
\)
Here \(\operatorname{Pf}(4,7)\) is the locus of \(7\times7\) skew-symmetric matrices
of rank at most \(4\), defined by the seven \(6\times6\) Pfaffians. We take a
complete intersection of \(\operatorname{Pf}(4,7)\) with thirteen general linear
sections of degree \(1\) to get a Fano 4-fold
\(
X\subset \mathbb P^7.
\)

\subsection{$\chi(T_X)$ of the Pfaffian subvarieties}

For the Pfaffian examples, we use the function
\codebreak{skewDegeneracyLocus}. Given a skew-symmetric morphism
\(\varphi:E\to E^\vee\otimes L\) on the ambient weighted projective space
\(P\), this function constructs the corresponding relative-Grassmannian model
\(Z\) of the degeneracy locus. In all examples,
the degeneracy locus has the expected codimension and contains no lower-rank
stratum, so \(Z\) is isomorphic to the desired Pfaffian subvariety \(X\). In Examples \(\#1,\#2,\#3,\#4\), and \(\#A\), we use
\codebreak{tangentBundle} to compute \(T_Z\), and apply \codebreak{chi} to compute
\(\chi(Z,T_Z)=\chi(X,T_X)\).

For Example \(\#6\), after constructing the Pfaffian degeneracy locus \(Z\)
using \codebreak{skewDegeneracyLocus}, we pull back the extra quartic line bundle
from the original weighted projective space to \(Z\). We then use
\codebreak{sectionZeroLocus} to impose this extra equation, and compute
\(\chi(X,T_X)\) by applying \codebreak{tangentBundle} and \codebreak{chi} to the
resulting zero locus.

The resulting values are as follows.
\[
\begin{array}{c|c}

\# & \chi(X,T_X) \\
\hline
1 & -80\\
2 & -81\\
3 & -94\\
4 & -77\\
6 & -181\\
A & -56\\

\end{array}
\]
\subsection{$\chi(T_X)$ of the weighted complete intersections}
\leavevmode\\[0.4em]
There are 13 families of smooth Fano fourfolds of index \(1\) and Picard
rank \(1\) realized as complete intersections in weighted projective spaces
\cite{kuchle,Przyjalkowski}. We write
\(X=X_{d_1,\ldots,d_c}\subset \mathbb P(a_0,\ldots,a_N)\). In each case,
\(X\) is smooth and does not meet the singular locus of the ambient weighted
projective space. We computed \(\chi(X,T_X)\) in Macaulay2 using
\codebreak{Schubert2}, together with the auxiliary function
\codebreak{weightedProjective\allowbreak Space}. Namely, after constructing
\(P=\mathbb P(a_0,\ldots,a_N)\), we define \(X\) as the zero locus of a
general section of
\(
\mathcal O_P(d_1)\oplus\cdots\oplus \mathcal O_P(d_c).
\)
Then \codebreak{sectionZeroLocus} gives \(X\), \codebreak{tangentBundle} gives
\(T_X\), and \codebreak{chi} computes \(\chi(X,T_X)\).

\[
\begin{array}{c|c|c|c}
\text{\#} & \mathbb P & \text{Degrees} & \chi(X,T_X) \\ \hline
1 & \mathbb P^6(1^3,2^2,3^2) & 6,6 & -166 \\
2 & \mathbb P^5(1^4,2,5) & 10 & -554 \\
3 & \mathbb P^6(1^4,2^2,3) & 4,6 & -197 \\
4 & \mathbb P^5(1^5,4) & 8 & -470 \\
5 & \mathbb P^5(1^5,2) & 6 & -255 \\
6 & \mathbb P^6(1^5,2^2) & 4,4 & -144 \\
7 & \mathbb P^6(1^6,3) & 2,6 & -320 \\
8 & \mathbb P^5 & 5 & -216 \\
9 & \mathbb P^6(1^6,2) & 3,4 & -145 \\
10 & \mathbb P^6 & 2,4 & -160 \\
11 & \mathbb P^6 & 3,3 & -116 \\
12 & \mathbb P^7 & 2,2,3 & -108 \\
13 & \mathbb P^8 & 2,2,2,2 & -84
\end{array}
\]
\vspace{1cm}
\subsection{$\chi(T_X)$ of the zero loci of homogeneous bundles on Grassmannians}
\leavevmode\\[0.4em]
There are 13 families of smooth Fano fourfolds with Picard rank 1, index 1 realized as zero loci of homogeneous vector bundles on Grassmannians \cite{kuchle}. See a list of 637 deformation families as zero loci of homogeneous vector bundles on products of Grassmannians by Bernardara, Fatighenti, Manivel, and Tanturri \cite{FanoFourfoldsDatabase} and its further discussion on K3 type \cite{fa}. In particular, the 13 families are the Picard rank 1 and index 1 cases in the database, and $\chi(X, T_X)$ has been calculated in the database. For completeness, we include the computation below. 

Let \(G=\operatorname{Gr}(a,n)\) denote the Grassmannian of
\(a\)-dimensional quotient spaces of \(\mathbb C^n\) and let $\mathcal U$ and $\mathcal Q$ denote the tautological subbundle and quotient bundle.  Each variety in the list is the zero locus of a general section of a homogeneous vector bundle $F$ constructed from $\mathcal U$ and $\mathcal Q$ on $G$.
Let $X=Z(s)\subset G$ be the zero locus of a general section of a homogeneous vector bundle $F$ on a Grassmannian $G$.  The calculations were implemented in Macaulay2 using the package
\codebreak{Schubert2}. For each family, we used the command \codebreak{sectionZeroLocus} to
construct \(X\) as the zero locus of a general section of \(\mathcal F\), then
used \codebreak{tangentBundle} to construct \(T_X\), and finally used
\codebreak{chi} to compute \(\chi(X,T_X)\).

\[
\begin{array}{c|c|c|c}
\text{K\"uchle type} & Y & \mathcal F & \chi(X,T_X) \\ \hline
(b1)  & \operatorname{Gr}(2,5) & \mathcal O_Y(3)\oplus \mathcal O_Y(1) & -109 \\
(b2)  & \operatorname{Gr}(2,5) & \mathcal O_Y(2)^{\oplus 2} & -72 \\
(b6)  & \operatorname{Gr}(2,6) & \mathcal O_Y(2)\oplus \mathcal O_Y(1)^{\oplus 3} & -63 \\
(c2)  & \operatorname{Gr}(3,6) & \wedge^2 \mathcal Q\oplus \mathcal O_Y(1)\oplus \mathcal O_Y(2) & -61 \\
(b10) & \operatorname{Gr}(2,7) & \mathcal U(1)\oplus \mathcal O_Y(2) & -62 \\
(b5)  & \operatorname{Gr}(2,6) & \mathcal Q(1)\oplus \mathcal O_Y(1)^{\oplus 2} & -48 \\
(b3)  & \operatorname{Gr}(2,6) & \wedge^3 \mathcal U(2) & -42 \\
(c1)  & \operatorname{Gr}(3,6) & \mathcal O_Y(1)^{\oplus 5} & -40 \\
(b8)  & \operatorname{Gr}(2,7) & \Sym^2 \mathcal Q\oplus \mathcal O_Y(1)^{\oplus 3} & -33 \\
(b11) & \operatorname{Gr}(2,8) & \mathcal U(1)\oplus \mathcal O_Y(1)^{\oplus 2} & -28 \\
(c6)  & \operatorname{Gr}(3,7) & \left(\wedge^2 \mathcal Q\right)^{\oplus 2}\oplus \mathcal O_Y(1)^{\oplus 2} & -28 \\
(c5)  & \operatorname{Gr}(3,7) & \wedge^2 \mathcal Q\oplus \mathcal U(1)\oplus \mathcal O_Y(1) & -25 \\
(c3)  & \operatorname{Gr}(3,7) & \mathcal U(1)^{\oplus 2} & -18
\end{array}
\]

\begin{remark}
    The Macaulay2 code for the computations of the Pfaffian subvarieties for the relative-Grassmannian method is included as an auxiliary file with the arXiv submission.
\end{remark}

\bibliographystyle{plain}
\bibliography{reference}

\end{document}